\documentclass[11pt,twoside,a4paper]{article}
\usepackage{amsmath,amssymb,amsthm}

\usepackage{graphicx,psfrag}
\newtheorem{thm}{Theorem}[section]
\newtheorem{lem}[thm]{Lemma}
\newtheorem{pro}[thm]{Proposition}

\theoremstyle{definition}

\theoremstyle{remark}

\newcommand{\R}{\mathbb{R}}

\newcommand{\N}{\mathbb{N}}

\newcommand{\cU}{\mathcal{U}}
\newcommand{\cV}{\mathcal{V}}

\newcommand{\al}{\alpha}

\newcommand{\ga}{\gamma}

\newcommand{\de}{\delta}

\newcommand{\si}{\sigma}

\newcommand{\la}{\lambda}
\newcommand{\La}{\Lambda}
\renewcommand{\phi}{\varphi}

\newcommand{\dist}{\operatorname{dist}}
\newcommand{\diam}{\operatorname{diam}}

\newcommand{\hyp}{\operatorname{H}}

\newcommand{\Int}{\operatorname{Int}}

\newcommand{\asdim}{\operatorname{asdim}}

\newcommand{\cone}{\operatorname{Co}}
\newcommand{\pt}{\operatorname{pt}}

\newcommand{\mesh}{\operatorname{mesh}}
\newcommand{\cdim}{\operatorname{cdim}}

\newcommand{\cp}{\operatorname{cap}}

\newcommand{\es}{\emptyset}
\renewcommand{\d}{\partial}
\newcommand{\di}{\d_{\infty}}
\newcommand{\set}[2]{\{#1:\,\text{#2}\}}
\newcommand{\sm}{\setminus}
\newcommand{\sub}{\subset}
\newcommand{\sups}{\supset}

\newcommand{\ov}{\overline}

\newcommand{\wh}{\widehat}

\begin{document}

\title{Capacity dimension and embedding of hyperbolic spaces
into the product of trees}
\author{Sergei Buyalo\footnote{Supported by RFFI Grant
05-01-00939 and Grant NSH-1914.2003.1}}

\date{}
\maketitle

\begin{abstract} We prove that every visual Gromov
hyperbolic space
$X$
whose boundary at infinity has the finite capacity
dimension,
$\cdim(\di X)<\infty$,
admits a quasi-isometric embedding into
$n$-fold
product of metric trees with
$n=\cdim(\di X)+1$.
\end{abstract}

\section{Introduction}\label{sect:introduct}

Recall that a map
$f:X\to Y$
between metric spaces is quasi-isometric if
$$\frac{1}{\La}|xx'|-\si\le|f(x)f(x')|\le\La|xx'|+\si$$
for some constants
$\La\ge 1$, $\si\ge 0$
and all
$x$, $x'\in X$.
Our main result is the following embedding theorem.

\begin{thm}\label{thm:main} Let
$X$
be a visual Gromov hyperbolic space, whose
boundary at infinity has the finite capacity
dimension,
$\cdim(\di X)<\infty$.
Then there exists a quasi-isometric embedding
$f:X\to T_1\times\dots\times T_n$
of
$X$
into
$n$-fold
product of metric trees
$T_1,\dots,T_n$
with
$n=\cdim(\di X)+1$.
\end{thm}

The property of a Gromov hyperbolic space
$X$
to be visual is a rough version of the property that every
point
$x\in X$
lies on a geodesic ray emanating from a fixed point
$x_0\in X$;
for the precise definition see sect.~\ref{subsect:visual}.

The boundary at infinity
$\di X$
is taken with a visual metric, which is defined up to
a quasi-symmetry of
$\di X$.
The notion of the capacity dimension of a metric space
is introduced in \cite{Bu}, where it is proved
that the capacity dimension is a quasi-symmetry invariant
(for a close notion of the Assouad-Nagata dimension this is
earlier proved in \cite{LS}).
In particular,
$\cdim(\di X)$
is independent of the choice of a visual metric on
$\di X$.

It is well-known that the asymptotic dimension of the
$n$-fold
product of any metric trees is at most
$n$
for every
$n\ge 1$
(see e.g. \cite[Chapter~9]{Ro}). Thus Theorem~\ref{thm:main}
generalizes the main result of \cite{Bu} saying that
$\asdim X\le\cdim\di X+1$
for every visual Gromov hyperbolic space
$X$.
In fact, Proposition~\ref{pro:charseqcov} below, which is the main
ingredient of the proof of Theorem~\ref{thm:main}, is
a refined version of \cite[Proposition~4.4]{Bu}.

In \cite{BS1}, it is constructed a quasi-isometric embedding
of the real hyperbolic space
$\hyp^n$
into the
$n$-fold
product of metric trees for every
$n\ge 2$.
It is known that
$\cdim\di\hyp^n=n-1$,
see \cite{Bu}, thus Theorem~\ref{thm:main} also generalizes
that result. Actually, our construction of the embedding is
a version of the construction from \cite{BS1}.

The estimate of the number of tree-factors, needed for a
quasi-isometric embedding given by Theorem~\ref{thm:main},
is sharp: according to \cite{BS2}, for every
$n\ge 2$,
the space
$\hyp^n$
admits no quasi-isometric embedding into the
$(n-1)$-fold
product of any metric trees stabilized by any Euclidean factor
$\R^N$.

The class of Gromov hyperbolic spaces to which Theorem~\ref{thm:main}
can be applied contains all visual Gromov
hyperbolic spaces with doubling boundary at infinity, or
which is the same, the spaces satisfying the bounded growth
condition at some scale, see \cite{BoS}, in particular, any
Gromov hyperbolic group and any Hadamard manifold with pinched
negative curvature is in this class. Indeed, the boundary
at infinity
$\di X$
of any such space can be quasi-symmetrically embedded in
$\R^N$
for some
$N\in\N$
by the Assouad embedding theorem \cite{As},
\cite[Chapter~12]{He}. Thus
$\cdim\di X\le\cdim\R^N=N$
is finite. The open problem is to find conditions, which
ensure that the capacity dimension of
$\di X$,
or more generally of a metric space, coincides with
the topological dimension.

For other similar embedding results into the product of
trees see \cite{Dr}, \cite{DZ}, \cite{LS} and references therein.

We briefly describe the structure of the paper.
Sect.~\ref{sect:capdim} is dedicated to the capacity dimension.
We fix notations and recall notions related to coverings of a
given metric space
$Z$
needed for the paper (sect.~\ref{subsect:cov}). Further,
in sect.~\ref{subsect:defcapdim}, we recall one of a number
of equivalent definitions of the capacity dimension, which
is suitable for our purposes.

The main feature of the capacity dimension is that the coverings
involved in its definition have the Lebesgue number at the
same scale as their mesh. This is the source of astonishing
flexibility in manipulating with coverings, which allows
to achieve many useful properties. The core of the paper
is sect.~\ref{subsect:charseq}, where we introduce the notion
of a
$\ga$-separated
characteristic sequence of coverings and prove the existence
of such a sequence under the condition
$\cdim Z<\infty$
(Proposition~\ref{pro:charseqcov}).

In sect.~\ref{sect:embd}, using a
$\ga$-separated
characteristic sequence of coverings of
$Z$,
we construct a quasi-isometric embedding of the hyperbolic
cone
$\cone(Z)$
over
$Z$
into the product of appropriate metric trees
(Theorem~\ref{thm:embcone}). The construction is a version
of the one from \cite{BS1} adopted to the current situation.

In sect.~\ref{sect:proofmain}, we recall basic notions from
the hyperbolic spaces theory and complete the proof of
Theorem~\ref{thm:main}.

\section{Capacity dimension}\label{sect:capdim}

Let
$Z$
be a metric space. For
$U$, $U'\sub Z$
we denote by
$\dist(U,U')$
the distance between
$U$
and
$U'$,
$\dist(U,U')=\inf\set{|uu'|}{$u\in U,\ u'\in U'$}$
where
$|uu'|$
is the distance between
$u$, $u'$.
For
$r>0$
we denote by
$B_r(U)$
the open
$r$-neighborhood
of
$U$, $B_r(U)=\set{z\in Z}{$\dist(z,U)<r$}$,
and by
$\ov B_r(U)$
the closed
$r$-neighborhood
of
$U$, $\ov B_r(U)=\set{z\in Z}{$\dist(z,U)\le r$}$.
We extend these notations over all real
$r$
putting
$B_r(U)=U$
for
$r=0$,
and defining
$B_r(U)$
for
$r<0$
as the complement of the closed
$|r|$-neighborhood
of
$Z\sm U$,
$B_r(U)=Z\sm\ov B_{|r|}(Z\sm U)$.
It is straightforward to check that the following is true.

\begin{lem}\label{lem:inclusion} Given
$U\sub Z$,
for every
$0<s<t$
we have
$$B_{t-s}(U)\sub B_{-s}(B_t(U)).$$
\qed
\end{lem}

A subset
$X$
of a metric space
$Z$
is a {\em net} in
$Z$
if there is
$\la>0$
such that
$\dist(z,X)\le\la$
for every
$z\in Z$.
In this case, we say that
$X$
is a
$\la$-net.

\subsection{Coverings}\label{subsect:cov}

Given a family
$\cU$
of subsets in a metric space
$Z$
we define
$\mesh(\cU)=\sup\set{\diam U}{$U\in\cU$}$.
The {\em multiplicity} of
$\cU$, $m(\cU)$,
is the maximal number of members of
$\cU$
with nonempty intersection. For
$r>0$,
the
$r$-{\em multiplicity}
of
$\cU$, $m_r(\cU)$,
is the multiplicity of the family
$\cU_r$
obtained by taking open
$r$-neighborhoods
of the members of
$\cU$.
So
$m_r(\cU)=m(\cU_r)$.
We say that a family
$\cU$
is {\em disjoint} if
$m(\cU)=1$,
and
$\cU$
is
$r$-{\em disjoint} if
$m_r(\cU)=1$.

A family
$\cU$
is called a {\em covering} of
$Z$
if
$\cup\set{U}{$U\in\cU$}=Z$.
A covering
$\cU$
is said to be {\em colored} if it is the union
of
$m\ge 1$
disjoint families,
$\cU=\cup_{a\in A}\cU^a$, $|A|=m$.
In this case we also say that
$\cU$
is
$m$-colored.
Clearly, the multiplicity of a
$m$-colored
covering is at most
$m$.

Let
$\cU$
be an open covering of a metric space
$Z$.
Given
$z\in Z$,
we let
$L'(\cU,z)=\sup\set{\dist(z,Z\sm U)}{$U\in\cU$}$,
$$L(\cU,z)=\min\{L'(\cU,z),\mesh(\cU)\}$$
be the Lebesgue number of
$\cU$
at
$z$
(the auxiliary
$L'(\cU,z)$
might be larger than
$\mesh(\cU)$
and even infinite as e.g. in the case
$Z=U$
for some member
$U\in\cU$),
$L(\cU)=\inf_{z\in Z}L(\cU,z)$
be the Lebesgue number of
$\cU$.
For this definition of the Lebesgue number, one needs
open coverings. For e.g. closed coverings, the Lebesgue
number is defined differently. We have
$L(\cU)\le L(\cU,z)\le\mesh(\cU)$
and for every
$z\in Z$
the open ball
$B_r(z)$
of radius
$r=L(\cU)$
centered at
$z$
is contained in some member of the covering
$\cU$.

We shall use the following obvious fact (see e.g. \cite{Bu}).

\begin{lem}\label{lem:insidecov} Let
$\cU$
be an open covering of
$Z$
with
$L(\cU)>0$.
Then for every
$s\in(0,L(\cU))$
the family
$\cU_{-s}=B_{-s}(\cU)$
is still an open covering of
$Z$,
and its
$s$-multiplicity
$m_s(\cU_{-s})\le m(\cU)$.
\qed
\end{lem}

\subsection{Definition of the capacity dimension}\label{subsect:defcapdim}

There are several equivalent definitions of the capacity
dimension, see \cite{Bu}. In this paper we shall use the
following one. Let
$\cU$
be an open covering of a metric space
$Z$.
We define the {\em capacity} of
$\cU$
by
$$\cp(\cU)=\frac{L(\cU)}
         {\mesh(\cU)}\in[0,1],$$
in the case
$\mesh(\cU)=0$
or
$L(\cU)=\mesh(\cU)=\infty$,
we put
$\cp(\cU)=1$
by definition. For
$\tau>0$, $\de\in(0,1)$
and an integer
$m\ge 0$
we put
$$c_{\tau}(Z,m,\de)=\sup_{\cU}\cp(\cU),$$
where the supremum is taken over all open,
$(m+1)$-colored
coverings
$\cU$
of
$Z$
with
$\de\tau\le\mesh(\cU)\le\tau$.

Next, we take
$$c(Z,m,\de)=\liminf_{\tau\to 0}c_{\tau}(Z,m,\de).$$
The function
$c(Z,m,\de)$
is monotone in
$\de$, $c(Z,m,\de')\ge c(Z,m,\de)$
for
$\de'<\de$.
Hence, there exists a limit
$c(Z,m)=\lim_{\de\to 0}c(Z,m,\de)$.
Now, we define the {\em capacity dimension} of
$Z$
as
$$\cdim(Z)=\inf\set{m}{$c(Z,m)>0$}.$$
In other words,
$\cdim(Z)$
is the minimal integer
$m\ge 1$
with the following property: There is a constant
$\de>0$
such that for every sufficiently small
$\tau>0$
there exists a
$(m+1)$-colored
open covering
$\cU$
of
$Z$
with
$\mesh(\cU)\le\tau$
and
$L(\cU)\ge\de\tau$.

\subsection{A characteristic sequence of coverings}\label{subsect:charseq}

Let
$Z$
be a metric space with finite capacity dimension,
$n=\cdim Z<\infty$.
We say that a sequence
$\cU_j$, $j\in\N$,
of
$(n+1)$-colored
(by a set
$A$)
open coverings of
$Z$
is {\em characteristic} with parameter
$r\in(0,1)$
if for some {\em characteristic} constants
$\de\in(0,1)$, $\la\ge 1$
the following conditions are fulfilled
\begin{itemize}
\item[(1)] $\mesh\cU_j\le r^j$
and
$L(\cU_j)\ge\de r^j$
for every
$j\in\N$;
\item[(2)] for every
$a\in A$, $j\in\N$,
the family
$\cU_j^a$
is a
$\la r^j$-net
in
$Z$.
\end{itemize}

Existence of a sequence of coverings with property (1)
follows directly from the definition
of the capacity dimension; property (2)
is auxiliary and easy to achieve, see Lemma~\ref{lem:charseqcov}.
A characteristic sequence of coverings
$\{\cU_j\}$
with parameter
$r$
is said to be {\em
$\ga$-separated}, $\ga\in(0,1)$,
if in addition the following holds

\begin{itemize}
\item[(3)] for every
$a\in A$
and for different members
$U\in\cU_j^a$, $U'\in\cU_{j\,'}^a$
with
$j\,'\le j$
we have either
$B_s(U)\cap U'=\es$,
or
$B_s(U)\sub U'$,
moreover, in the case
$j\,'<j$,
there is
$U''\in\cU_j^a$
with
$B_s(U'')\sub U'$
for
$s=\ga r^j$.
\end{itemize}

Sequences of coverings with variants of (1)--(3)
have been used in a number of papers, see e.g.
\cite{Bu}, \cite{BS1}, \cite{Dr}, \cite{LS}, for various
purposes, basically for constructing embeddings with
specific properties.

\begin{pro}\label{pro:charseqcov} Suppose that
$Z$
is a metric space with finite capacity dimension,
$\cdim Z\le n$.
Then there are constants
$\de$, $\ga\in(0,1)$, $\la\ge 1$
such that for every sufficiently small
$r>0$
there exists a
$\ga$-separated
characteristic sequence of coverings
$\cU_j$, $j\in\N$,
of
$Z$
with the parameter
$r$
and characteristic constants
$\de$, $\la$.
\end{pro}

The proof of this Proposition follows the same line as
\cite[Proposition~4.4]{Bu}, cf. also \cite[Proposition~4.1]{LS}.
We first construct a characteristic sequence of coverings of
$Z$
and then modify it to obtain property (3).

\begin{lem}\label{lem:charseqcov} Under the condition of
Proposition~\ref{pro:charseqcov}, there are constants
$\de\in(0,1)$, $\la\ge 1$
such that for every sufficiently small
$r>0$
there exists a characteristic sequence of coverings
$\wh\cU_j$, $j\in\N$,
of
$Z$
with the parameter
$r$
and characteristic constants
$\de$, $\la$.
Moreover, for every
$j\in\N$,
every
$U\in\wh\cU_j$
contains a ball of radius
$\de r^j$
and for every
$a\in A$
the family
$\wh\cU_j^a$
is
$\de r^j$-disjoint.
\end{lem}

\begin{proof}
We have
$c=c(Z,n)/8>0$
by the definition of
$\cdim Z$, $c(Z,n,\de')\ge 4c$
for all sufficiently small
$\de'>0$.
We fix such a
$\de'$,
then
$c_{\tau}(Z,n,\de')\ge 2c$
for all
$\tau$, $0<\tau\le\tau_0$.
This means that for every
$\tau\in(0,\tau_0]$
there is a
$(n+1)$-colored
open covering
$\cU_\tau$
of
$Z$
with
$\de'\tau\le\mesh(\cU_\tau)\le\tau$
and the capacity arbitrarily close to
$c_{\tau}(Z,n,\de')$,
in particular,
$L(\cU_\tau)\ge c\mesh(\cU_\tau)\ge c\de'\tau$.

Take a positive
$r<r_0=\min\{c\de'/2,\tau_0\}$
and for every
$j\in\N$
consider the covering
$\cU_j=\cU_{\tau_j}$,
where
$\tau_j=r^j$.
Then the sequence
$\cU_j$, $j\in\N$,
of
$(n+1)$-colored
coverings of
$Z$
satisfies condition (1) (with characteristic constant
$c\de'$).
Fix
$j\in\N$
and for
$s=c\de'r^j/2$
consider the family
$\wh\cU_j=B_{-s}(\cU_j)$.
By Lemma~\ref{lem:insidecov},
$\wh\cU_j$
is an open covering of
$Z$,
and we have
$\mesh(\wh\cU_j)\le\mesh(\cU_j)\le r^j$,
$L(\wh\cU_j)\ge\frac{1}{2}c\de'r^j$,
thus (1) is satisfied for
$\wh\cU_j$
with
$\de=c\de'/2$.

We can additionally assume that every
$U\in\wh\cU_j$
contains a ball of radius
$\de r^j\le L(\wh\cU_j)$
for otherwise
$U$
is covered by other members of
$\wh\cU_j$
and thus it can be deleted from
$\wh\cU_j$
without destroying property (1).

Furthermore, for every color
$a\in A$,
the family
$\wh\cU_j^a$
is
$\de r^j$-disjoint.
Now, adding to a given family
$\wh\cU_j^a$, $a\in A$,
copies of members of other colors from
$\wh\cU_j$,
which are
$\de r^j$-disjoint
with
$\wh\cU_j^a$,
we can make out of it a
$\la r^j$-net
with
$\la=1+2\de$.
This does not change the
$\mesh$
and does not decrease the Lebesgue number of
$\wh\cU_j$.
Thus (2) is satisfied for every
$\wh\cU_j$, $j\in\N$.
\end{proof}

We shall use the following modification of
a construction from \cite{Bu}.
Let
$\cU$, $\cU'$
be families of sets in
$Z$, $s>0$.
We denote by
$\cU\ast_s\cU'$
the family obtained by taking for every
$U\in\cU$
the
$s$-neighborhood
of the union
$V$
of
$U$
and all members
$U'\in\cU'$
with nonempty
$B_s(U)\cap B_s(U')$,
$\cU\ast_s\cU'=\set{B_s(V)}{$U\in\cU$}$.

\begin{lem}\label{lem:undisfam} Assume that a family of sets
$\wh\cU$
is
$\de s$-disjoint
for
$\de\in(0,2/3]$, $s>0$,
and
$\mesh(\wh\cU)\le 2s$.
Then, the operation
$$U\mapsto U^\ast=B_{-4s}(U)\ast_{\de s}\wh\cU$$
does not increase any set
$U\sub Z$, $U^\ast\sub U$,
and for every
$\wh U\in\wh\cU$,
it holds either
$B_{\de s}(\wh U)\cap U^\ast=\es$
or
$B_{\de s}(\wh U)\sub U^\ast$.
\end{lem}

\begin{proof}
For every
$\wh U\in\wh\cU$,
we have
$\diam B_{\de s}(\wh U)\le 2s+2\de s$,
thus
$\de s+\diam B_{\de s}(\wh U)\le 2s(1+3\de/2)\le 4s$.
Hence if
$B_{\de s}(\wh U)$
intersects
$B_{\de s}(B_{-4s}(U))$
then
$B_{\de s}(B_{-4s}(U)\cup\wh U)\sub U$,
and
$U^\ast\sub U$,
in particular, in that case
$B_{\de s}(\wh U)\sub U^\ast$.
Otherwise,
$B_{\de s}(\wh U)$
misses the
$\de s$-neighborhood
of
$B_{-4s}(U)$
as well as the one of every other member of
$\wh\cU$,
and thus
$B_{\de s}(\wh U)\cap U^\ast=\es$.
\end{proof}

\begin{lem}\label{lem:undisfam1} Under the conditions of
Lemma~\ref{lem:undisfam}, assume that
$B_t(U_1)\sub U_2$, $t>4s$,
for some sets
$U_1$, $U_2\sub Z$.
Then
$B_{t'}(U_1^\ast)\sub U_2^\ast$
for
$t'=t-4s$.
\end{lem}

\begin{proof} Using Lemma~\ref{lem:undisfam} we obtain
$B_{t'}(U_1^\ast)\sub B_{t'}(U_1)$.
By Lemma~\ref{lem:inclusion},
$B_{t-4s}(U_1)\sub B_{-4s}(B_t(U_1))$.
Finally, we have
$B_{-4s}(B_t(U_1))\sub B_{-4s}(U_2)\sub U_2^\ast$.
\end{proof}

\begin{proof}[Proof of Proposition~\ref{pro:charseqcov}]
Using Lemma~\ref{lem:charseqcov}, we find a characteristic
sequen\-ce of coverings
$\wh\cU_j$, $j\in\N$
with characteristic constants
$\de\in(0,1)$, $\la\ge 1$
and arbitrarily small parameter
$r$
such that for every
$j\in\N$,
every
$U\in\wh\cU_j$
contains a ball of radius
$\de r^j$
and for every color
$a\in A$,
the family
$\wh\cU_j^a$
is
$\de r^j$-disjoint.
We further assume that
$\frac{2r}{1-r}\le\de/4\le 1/6$
and
$(\la+1)r<\de/2$.

Fix a color
$a\in A$
and define
$\cV_1^a=\cU_{1,1}^a:=\wh\cU_1^a$.
Then, the family
$\cU_{1,1}^a$
is
$\de r$-disjoint
and
$\mesh(\cU_{1,1}^a)\le r$.

Assume that for
$k\ge 1$
the family
$\cV_k^a$
is already defined, and it has the following properties
\begin{itemize}
\item[(i)] $\cV_k^a=\cup_{j=1}^k\cU_{j,k}^a$;
\item[(ii)] for every
$1\le j\le k$,
the family
$\cU_{j,k}^a$
is
$\de r^j$-disjoint
and
$\mesh(\cU_{j,k}^a)\le r^j$;
\item[(iii)] given
$1\le j\,'<j\le k$,
for every
$U'\in\cU_{j\,',k}^a$, $U\in\cU_{j,k}^a$,
we have either
$B_t(U)\cap U'=\es$
with
$t=\de r^j/2$
or
$B_t(U)\sub U'$
with
$t=\ga_{k,j}r^j$,
where
$\ga_{k,j}$
is defined recurrently by
$\ga_{j,j}=\de/2$
and
$\ga_{k,j}=\ga_{k-1,j}-2r^{k-j}$
for
$k>j$.
\end{itemize}
We define
$$\cV_{k+1}^a:=B_{-4s}(\cV_k^a)
  \ast_{\de s}\wh\cU_{k+1}^a\cup\wh\cU_{k+1}^a$$
with
$s=r^{k+1}/2$.
Then
$\cV_{k+1}^a=\cup_{j=1}^{k+1}\cU_{j,k+1}^a$,
where
$$\cU_{j,k+1}^a=B_{-4s}(\cU_{j,k}^a)\ast_{\de s}\wh\cU_{k+1}^a$$
for
$1\le j\le k$
and
$\cU_{k+1,k+1}^a=\wh\cU_{k+1}^a$.
Since the family
$\wh\cU_{k+1}^a$
is
$2\de s$-disjoint,
we can apply Lemma~\ref{lem:undisfam}, by which every
$U^\ast\in\cU_{j,k+1}^a$
with
$1\le j\le k$
is contained in the appropriate
$U\in\cU_{j,k}^a$,
in particular, the family
$\cU_{j,k+1}^a$
is
$\de r^j$-disjoint
and
$\mesh(\cU_{j,k+1}^a)\le\mesh(\cU_{j,k}^a)\le r^j$.
Furthermore, for every
$\wh U\in\cU_{k+1,k+1}^a$
we have either
$B_{\de s}(\wh U)\cap U^\ast=\es$
or
$B_{\de s}(\wh U)\sub U^\ast$.

Now, if
$U'\in\cU_{j\,',k}^a$
with
$j\,'<j$
and
$B_t(U)\cap U'=\es$
with
$t=\de r^j/2$
then
$B_t(U^\ast)\cap U'^\ast=\es$
since
$U^\ast\sub U$
and
$U'^\ast\sub U'$.
In the case
$B_t(U)\sub U'$
with
$t=\ga_{j,k}r^j$
by Lemma~\ref{lem:undisfam1} we have
$B_{t-2r^{k+1}}(U^\ast)\sub U'^\ast$
and
$t-2r^{k+1}=\ga_{k+1,j}r^j$
with
$\ga_{k+1,j}=\ga_{k,j}-2r^{k+1-j}$.
Note that
$\lim_{k\to\infty}\ga_{k,j}=\de/2-\frac{2r}{1-r}\ge\de/4$
for every
$j>1$.

Therefore, for every color
$a\in A$,
we have the sequence
$\cV_k^a$, $k\in\N$,
of families of sets in
$Z$
with properties (i)--(iii). It follows from the definition
of the
$\ast$-operation
that every member
$U^\ast\in\cV_{k+1}^a$
is contained in its well defined predecessor
$U\in\cV_k^a$,
moreover,
$U^\ast\in\cU_{j,k+1}^a$
if and only if
$U\in\cU_{j,k}^a$.
In this sense, the sequence
$\cV_k^a$
is monotone,
$\cV_k^a\sups\cV_{k+1}^a$,
and we define
$\cU_j^a=\Int(\cap_{k\ge j}\cU_{j,k}^a)$,
$\cU_j=\cup_{a\in A}\cU_j^a$
for every
$j\in\N$.

We put
$\wh s_j=\sum_{k\ge j}2r^{k+1}=\frac{2r}{1-r}r^j$.
Then
$\wh s_j\le\de r^j/4<L(\wh\cU_j)$
and
$B_{-\wh s_j}(\wh\cU_j^a)\sub \cU_j^a$
for every
$a\in A$, $j\in\N$.
By Lemma~\ref{lem:insidecov}, the family
$\cU_j$
is still an open
$(n+1)$-colored
covering of
$Z$
with
$L(\cU_j)\ge \de r^j-\wh s_j\ge\de r^j/2$.
From (ii), we obtain
$\mesh(\cU_j)\le r^j$
and the family
$\cU_j^a$
is
$\de r^j$-disjoint
for every
$a\in A$, $j\in\N$.
Property (iii) implies that given
$1\le j\,'<j$,
for every
$U'\in\cU_{j\,'}^a$, $U\in\cU_j^a$
we have either
$B_t(U)\cap U'=\es$
or
$B_t(U)\sub U'$
with
$t=\ga r^j$, $\ga\ge\de/4$.

Finally, recall that for every
$j\in\N$, $a\in A$,
the family
$\wh\cU_j^a$
is a
$\la r^j$-net
with
$\la\ge 1$
independent of
$j$,
and every
$\wh U\in\wh\cU_j$
contains a ball of radius
$\de r^j$.
Since every
$U\in\cU_j$
contains some
$\wh U\in B_{-\wh s_j}(\wh\cU_j)$,
it also contains a ball
$B\sub \wh U$
of radius
$\de r^j-\wh s_j\ge\de r^j/2$.
Therefore, for every color
$a\in A$
and every
$z\in Z$,
there is
$U\in\cU_j^a$
with
$\dist(z,U)\le(\la+1)r^j$,
which means that the family
$\cU_j^a$
is
$(\la+1)r^j$-net.
Since
$(\la+1)r<\de/2$,
this also shows that for every
$U'\in\cU_{j\,'}^a$
with
$j\,'<j$
there is
$U''\in\cU_j^a$
with
$U''\cap U'\neq\es$,
hence
$B_{\ga r^j}(U'')\sub U'$.
This completes the proof of Proposition~\ref{pro:charseqcov}.
\end{proof}

\section{The embedding construction}\label{sect:embd}

\subsection{The hyperbolic cone}\label{subsect:hypcone}

Let
$Z$
be a bounded metric space. Assuming that
$\diam Z>0$,
we put
$\mu=\pi/\diam Z$
and note that
$\mu|zz'|\in[0,\pi]$
for every
$z$, $z'\in Z$.
Recall that the hyperbolic cone
$\cone(Z)$
over
$Z$
is the space
$Z\times[0,\infty)/Z\times\{0\}$
with metric defined as follows. Given
$x=(z,t)$, $x'=(z',t')\in\cone(Z)$
we consider a triangle
$\ov o\,\ov x\,\ov x'\sub\hyp^2$
with
$|\ov o\,\ov x|=t$, $|\ov o\,\ov x'|=t'$
and the angle
$\angle_{\ov o}(\ov x,\ov x')=\mu|zz'|$.
Now, we put
$|xx'|:=|\ov x\,\ov x'|$.
In the degenerate case
$Z=\{\pt\}$,
we define
$\cone(Z)=\{\pt\}\times[0,\infty)$
as the metric product. The point
$o=Z\times\{0\}\in\cone(Z)$
is called the {\em vertex} of
$\cone(Z)$.

\begin{thm}\label{thm:embcone} Let
$Z$
be a bounded metric space with finite
capacity dimension,
$\cdim Z=n<\infty$.
Then there exists a quasi-isometric embedding
$f:\cone(Z)\to\prod_{a\in A}T_a$
of
$\cone(Z)$
into the
$(n+1)$-fold
product of metric trees
$T_a$, $|A|=n+1$.
\end{thm}

The construction of an embedding
$f:\cone(Z)\to\prod_{a\in A}T_a$,
we are using in this paper, is similar to the one from \cite{BS1}.

\subsection{Construction of trees
$T_a$}\label{subsect:trees}
Every tree
$T_a$, $a\in A$
is a rooted simplicial tree with a root
$v_a\in T_a$,
and every edge of
$T_a$
has length 1.
Using Proposition~\ref{pro:charseqcov}, we find
constants
$\de$, $\ga\in(0,1)$, $\la\ge 1$,
and, for a sufficiently small
$r>0$,
a
$\ga$-separated
characteristic sequence of open
$(n+1)$-colored
coverings
$\cU_j=\cup_{a\in A}\cU_j^a$, $j\in\N$,
of
$Z$
with the parameter
$r$
and characteristic constants
$\de$, $\la$.
We assume that
$\la r<\de$
and
$r<\diam Z$.

For every
$a\in A$,
we define a graph
$T_a$
as follows. Its vertex set
$V^a$
is the disjoint union,
$V^a=\cup_{j\ge 0}V_j^a$,
where the vertex set
$V_j^a$
of level
$j\ge 0$
is identified with the set
$\cU_j^a$, $\cU_0^a=\{Z\}$,
i.e. the root
$v_a$
corresponds to the set
$Z$.
There are edges between only vertices of distinct levels,
and vertices
$v\in V_j^a$, $v'\in V_{j\,'}^a$, $j\,'<j$,
are connected by the (unique) edge if and only if
for the corresponding members
$U\in\cU_j^a$, $U'\in\cU_{j\,'}^a$
one holds
$U\sub U'$
and
$j\,'$
is maximal with this property.

By this definition, every vertex
$v\in V_j^a$, $j\in\N$,
is connected with a lower lever vertex
$v'$
by an edge, and every vertex from
$V_1^a$
is connected with
$v_a$.
Thus, the graph
$T_a$
is connected. By properties of the
sequence
$\{\cU_j\}$,
for every vertex
$v\in V_j^a$, $j\in\N$,
there is at most one edge leading to a lower level
vertex, therefore,
$T_a$
is a tree.

\subsection{Construction of an embedding
$f:\cone(Z)\to\prod_{a\in A}T_a$}\label{subsect:emb}
We denote by
$Z_t$
the metric sphere of radius
$t>0$
around
$o$
in
$\cone(Z)$.
There are natural polar coordinates
$x=(z,t)$, $z\in Z$, $t\ge 0$,
in
$\cone(Z)$.
Then
$Z_t=\set{(z,t)}{$z\in Z$}$
is the copy of
$Z$
at the level
$t$.
For
$t>0$
we denote by
$\pi_t:Z_t\to Z$
the canonical homeomorphism,
$\pi_t(z,t)=z$.

We denote
$R=\ln\frac{1}{r}$
and for
$j\in\N$
put
$Z_j=Z_{jR}$,
$\pi_j=\pi_{jR}$.
We let
$Z_0=o$
be the vertex of
$\cone(Z)$
and consider the set
$X=\cup_{j\ge 0}Z_j\sub\cone(Z)$.
Given
$a\in A$,
we define
$f_a:X\to T_a$
as follows. We put
$f_a(o)=v_a$
and for
$x\in Z_j$, $j\in\N$,
we take a member
$U\in\cU_j^a$
closest to
$\pi_j(x)\in Z$,
and let
$f_a(x)=v\in V_j^a$
be the vertex corresponding to
$U$.

We need to pass from distances in
$Z$
to distances in any sphere
$Z_j$, $j\in\N$.
This is done in important technical
Lemma~\ref{lem:projconedist}. We say that
$A\asymp B$
up to a multiplicative error
$\le C$,
if
$$\frac{1}{C}\le\frac{A}{B}\le C.$$

\begin{lem}\label{lem:projconedist} Given
$z$, $z'\in Z$, $j\in\N$
we let
$\mu|zz'|=\tau$,
$\tau_j=|z_jz_j'|$,
where
$z_j=\pi_j^{-1}(z)$, $z_j'=\pi_j^{-1}(z')\in Z_j$.
Then
$$\sinh(\tau_j/2)\asymp\tau/r^j$$
up to a universally bounded multiplicative error.
\end{lem}

\begin{proof}
Using the hyperbolic cosine law
$\cosh a(t,\al)=\cosh^2(t)-\sinh^2(t)\cos\al$
for the base
$a(t,\al)$
of an isosceles triangle in
$\hyp^2$
with sides
$t$
and the angle
$\al$
between them, we see that
\begin{eqnarray*}
 \cosh\tau_j&=&\cosh^2(jR)-\sinh^2(jR)
 \cos\tau\\
  &=&1+\sinh^2(jR)2\sin^2\left(\frac{\tau}{2}\right)
\end{eqnarray*}
by the definition of the cone metric. The claim follows.
\end{proof}

\subsubsection{The large scale Lipschitz property of
$f_a$}\label{subsubsect:largelip}

The proof of the large scale Lipschitz property of
$f_a$
follows the same line as \cite[Proposition~2.7]{BS1}.
We need the following elementary fact, see \cite[Lemma~2.6]{BS1}.

\begin{lem}\label{lem:trisimple} Let
$p$, $q$, $t$
be the side lengths of a triangle in a metric space
such that
$t\ge p$.
Then
$p+q\le 3t$.
\qed
\end{lem}

\begin{pro}\label{pro:roughlip} For every
$a\in A$,
the map
$f_a:X\to T_a$
is roughly Lipschitz, i.e. there are constants
$\La>0$, $\si\ge 0$,
such that
$$|f_a(x)f_a(x')|\le\La|xx'|+\si$$
for each
$x$, $x'\in X$.
\end{pro}

\begin{proof} Fix
$x$, $x'\in X$
and consider
$v=f_a(x)$, $v'=f_a(x')$.
We can assume without loss of generality that
$x\in Z_j$, $x'\in Z_{j\,'}$
with
$j\ge j\,'\ge 0$.
Then,
$v\in V_j^a$, $v'\in V_{j\,'}^a$,
and we also assume that
$v\neq v'$.

First, we consider the case
$j\,'=j$.
It follows from properties of edges that
no shortest path in
$T_a$
has an interior vertex with locally maximal level. Thus
the shortest path in
$T_a$
between
$v$
and
$v'$
has a unique vertex
$v_0$
of a lowest level
$j_0$, $j_0<j$.
On the geodesic segments
$v_0v$, $v_0v'\sub T_a$,
we take the vertices
$v_1\in v_0v$, $v_1'\in v_0v'$
adjacent to
$v_0$.
We can assume that
$v_1\in V_{j_1}^a$, $v_1'\in V_{j_1'}^a$,
where
$j_0<j_1'\le j_1\le j$.

For the covering member
$U_1\in\cU_{j_1}^a$, $U_1'\in\cU_{j_1'}^a$
corresponding to
$v_1$, $v_1'$
respectively, we have either
$B_s(U_1)\cap U_1'=\es$
or
$B_s(U_1)\sub U_1'$
for
$s=\ga r^{j_1}$
by the separation property (property (3))
of the sequence
$\{\cU_j\}$.
The last possibility is excluded since otherwise
there is a path in
$T_a$
between
$v$
and
$v'$,
missing
$v_0$,
and hence the initial path
$vv_0\cup v_0v'$
is not the shortest one, a contradiction.
Therefore,
$\dist(U_1,U_1')\ge\ga r^{j_1}$
in
$Z$,
and by Lemma~\ref{lem:projconedist}, we have
$$\dist(\wh U_1,\wh U_1')\ge2R(j-j_1)-\si'$$
for some constant
$\si'$
depending only on
$\mu\ga$,
where
$\wh U_1=\pi_j^{-1}(U_1)$, $\wh U_1'=\pi_j^{-1}(U_1')\sub Z_j$.
We have
$\dist(\pi_j(x),U_1)\le\dist(\pi_j(x),U)\le\la r^j$,
where
$U\in\cU_j^a$
corresponds to the vertex
$f_a(x)\in T_a$,
and similarly,
$\dist(\pi_j(x'),U_1')\le\la r^j$.
Thus using Lemma~\ref{lem:projconedist}, we obtain
$|xx'|\ge\dist(\wh U_1,\wh U_1')-2\la'$,
with
$\la'>0$
depending only on
$\la\mu$.

Let
$j_2'$
be the level of the vertex
$v_2'\in v_1'v'$
adjacent to
$v_1'$, $v_2'\in V_{j_2'}^a$.
There are two possibilities, which are treated differently.

(1) $j_2'\ge j_1$.
Then for the distances in the tree
$T_a$,
we have
$|v_0v|=|v_1v|+1\le j-j_1+1$
and
$|v_0v'|=|v_2'v'|+2\le j-j_2'+2\le j-j_1+2$.
Therefore,
$$|vv'|\le 2(j-j_1)+3\le\La|xx'|+\si$$
for
$\La=1/R$
and some
$\si\ge 0$
independent of
$x$, $x'$.

(2) $j_2'<j_1$.
Let
$U_2'\in\cU_{j_2'}^a$
be the member corresponding to the vertex
$v_2'$.
By the separation property,
$B_s(U_2')\sub U_1'$
for
$s=\ga r^{j_2'}$.
Since
$U_1\cap U_1'=\es$,
we have
$\dist(U_1,U_2')\ge\ga r^{j_2'}$.
Applying Lemma~\ref{lem:projconedist}, we obtain
$\dist(\wh U_1,\wh U_2')\ge 2R(j-j_2')-\si'$,
where
$\wh U_2'=\pi_j^{-1}(U_2')$.
Similarly to (1), this yields
$|vv'|\le\La|xx'|+\si$.

Finally, we consider the case
$j'<j$.
Let
$U'\in\cU_{j\,'}^a$
be the member corresponding to the vertex
$v'$, $\wh U'=\pi_{j\,'}^{-1}(U')\sub Z_{j\,'}$.
Since
$\dist(\pi_{j\,'}(x'),U')\le\la r^{j\,'}$
and
$\diam U'\le r^{j\,'}$,
we have
$|x'y|\le\la''$
for every
$y\in\wh U'$
by Lemma~\ref{lem:projconedist}, where the constant
$\la''$
depends only on
$\la$, $\mu$.
Thus without loss of generality, we can take as
$x'$
any point from
$\wh U'$.
Recall that by property (3) of the sequence
$\{\cU_j\}$,
there is
$U''\in\cU_j^a$
with
$B_s(U'')\sub U'$
for
$s=\ga r^j$.
We take
$x'\in\wh U'$
so that
$\pi_{j\,'}(x')\in U''$.

Now, consider
$x''\in Z_j$
sitting over
$x'$,
i.e.
$\pi_j(x'')=\pi_{j\,'}(x')$.
Then
$x''\in\wh U''=\pi_j^{-1}(U'')$
and thus
$f_a(x'')=v''\in V_j^a$
corresponds to
$U''$.
It follows that
$v'$
is the lowest level vertex of the segment
$v'v''\sub T_a$,
and
$|v'v''|\le j-j'=\La|x'x''|$.
By what we have already proved,
$|vv''|\le\La|xx''|+\si$.
On the other hand, obviously,
$|x'x''|\le|x'x|$,
and thus
$|xx''|+|x''x'|\le 3|xx'|$
by Lemma~\ref{lem:trisimple}. Therefore,
$$|vv'|\le|vv''|+|v''v'|\le\La(|xx''|+|x''x'|)+\si
  \le3\La|xx'|+\si,$$
which completes the proof of the Proposition.
\end{proof}

\subsubsection{The large scale bilipschitz property of
$f$}\label{subsubsect:largebilip}

The map
$f:X\to\prod_{a\in A}T_a$
defined by its coordinate maps
$f_a:X\to T_a$
is roughly Lipschitz by Proposition~\ref{pro:roughlip}.
To prove that
$f$
is roughly bilipschitz, we begin with the following Lemma,
which is the main ingredient of the proof.

\begin{lem}\label{lem:radial} Given
$x=(z,Rj)\in Z_j$, $j\in\N$,
for every nonnegative integer
$i\le j$,
there is a color
$a\in A$
such that
$\dist(f_a(x),V_i^a)\ge M$
with
$M+1\ge(j-i+1)/|A|$.
Furthermore, if for
$k\le i$
the vertex
$v\in V_k^a$
is the lowest level vertex of the segment
$f_a(x)v\sub T_a$,
then
$|f_a(x)v|\ge M$.
\end{lem}

\begin{proof} Recall that the Lebesgue number
$L(\cU_j)\ge\de r^j$
for every integer
$j\ge 0$,
since the sequence of coverings
$\{\cU_j\}$
is characteristic. Thus for every
$j$,
there is
$U_j\in\cU_j$
with
$B_{\de r^j}(z)\sub U_j$.
There is a color
$a\in A$
such that the set
$\{U_i,\dots,U_j\}$
contains
$M+1\ge(j-i+1)/|A|$
members having the color
$a$,
i.e. every of those
$U_k\in\cU_k^a$.
Let
$U\in\cU_j^a$
be the member corresponding to the vertex
$f_a(x)\in T_a$.
Since
$\la r<\de$,
we have
$\dist(z,U)\le\la r^j<\de r^k$
for every
$k<j$.
Thus
$B_{\de r^k}(z)\cap U\neq\es$
and, by the separation property,
$U\sub U_k$
for every
$k<j$
with
$U_k$
having the color
$a$.

Using again the separation property, we obtain
that any path in
$T_a$
between
$f_a(x)$
and the set
$V_i^a$
must contain at least
$M+1$
vertices and hence
$\dist(f_a(x),V_i^a)\ge M$.

Finally, let
$V\in\cU_k^a$
be the set corresponding to the vertex
$v$.
By the assumption on
$v$,
the set
$V$
contains
$U$
and every set from the list
$\{U_i,\dots,U_j\}$
having the color
$a$.
Hence,
$|f_a(x)v|\ge M$.
\end{proof}

\begin{pro}\label{pro:bilip} There are constants
$\La>0$, $\si\ge 0$
such that
$$|xx'|\le\La|f(x)f(x')|+\si$$
for all
$x$, $x'\in X$.
\end{pro}

\begin{proof} We have
$x=(z,Rj)$, $x'=(z',Rj\,')$
in polar coordinates in
$\cone(Z)$.
We can assume that
$j\ge j\,'$.
Furthermore, if
$j\,'=0$
then
$x'=o$,
the vertex of
$\cone(Z)$,
and we assume
$z'=z$
in this case.

First, consider the case
$|zz'|\le\la r^j+(\la+1)r^{j\,'}$.
For
$x''=(z,Rj\,')$
we have
$|x''x'|\le\si'$
by Lemma~\ref{lem:projconedist},
where
$\si'$
depends only on
$\la$
and
$\mu$.
Then
$|xx'|\le|xx''|+|x''x'|\le(j-j\,')R+\si'$.

On the other hand,
$f_a(x')\in V_{j\,'}^a$
for all
$a\in A$,
and by Lemma~\ref{lem:radial}, we have
$$|f_a(x)f_a(x')|\ge\dist(f_a(x),V_{j\,'}^a)\ge M$$
with
$M+1\ge(j-j\,'+1)/|A|$
for some color
$a\in A$.
Therefore,
$|xx'|\le\La|f(x)f(x')|+\si$
with
$\La=|A|R$
and
$\si$
depending only on
$|A|$, $R$, $\si'$.

Second, consider the case
$|zz'|>\la r^j+(\la+1)r^{j\,'}$.
There is an integer
$l\ge 0$
with
$$\la r^j+r^{l+1}+\la r^{j\,'}<|zz'|\le
  \la r^j+r^l+\la r^{j\,'},$$
where we can take
$r^0:=\diam Z$,
since
$r<\diam Z$.
Then
$r^{l+1}\ge r^{j\,'}$
and thus
$l+1\le j\,'$.
Now, we show that for every color
$a\in A$,
any path in
$T_a$
between
$f_a(x)$
and
$f_a(x')$
passes through a vertex of a level
$k\le l$.
Indeed, let
$v\in T_a$
be a lowest level vertex of the segment
$f_a(x)f_a(x')\sub T_a$, $v\in V_k^a$.
Then, the set
$V\in\cU_k^a$
corresponding to
$v$
contains both
$U\in\cU_j^a$
and
$U'\in\cU_{j\,'}^a$
corresponding to
$f_a(x)$, $f_a(x')$
respectively, and we have
$$|zz'|\le\dist(z,U)+\diam V+\dist(z',U')
  \le\la r^j+r^k+\la r^{j\,'}.$$
It follows
$r^{l+1}<r^k$
and thus
$k\le l$.

By Lemma~\ref{lem:radial}, there is a color
$a\in A$
such that
$\dist(f_a(x),V_l^a)\ge M$
with
$M+1\ge(j-l+1)/|A|$.
Let
$v\in V_k^a$
be the lowest level vertex of the segment
$f_a(x)f_a(x')$.
We have
$|f_a(x)f_a(x')|\ge|f_a(x)v|$,
and by Lemma~\ref{lem:radial},
$|f_a(x)v|\ge M$
since
$k\le l$.

Consider the points
$x_l=(z,Rl)$, $x_l'=(z',Rl)\in Z_l$.
Using the estimate
$$|zz'|\le\la r^j+r^l+\la r^{j\,'}\le(2\la+1)r^l$$
and Lemma~\ref{lem:projconedist}, we obtain
$|x_lx_l'|\le\si'$
with
$\si'$
depending only on
$\la$
and
$\mu$.
Thus
$$|xx'|\le|xx_l|+|x_lx_l'|+|x_l'x'|
  \le 2(j-l)R+\si'\le\La|f(x)f(x')|+\si,$$
where
$\La=2|A|R$,
and
$\si$
depends only on
$|A|$, $R$, $\si'$.
\end{proof}

\subsection{Proof of Theorem~\ref{thm:embcone}}\label{subsect:embcone}

The map
$f:X\to\prod_{a\in A}T_a$
is a quasi-isometric embedding by Propositions~\ref{pro:roughlip}
and \ref{pro:bilip}. The set
$X\sub\cone(Z)$
being the union of equidistant spheres
$Z_j$, $j\ge 0$,
is obviously
$R/2$-net
in
$\cone(Z)$, $R=\ln\frac{1}{r}$.
Thus to define a quasi-isometric embedding of
$\cone(Z)$
into something it suffices to define it on
$X$.
Therefore, the map
$f$
is a well defined quasi-isometric embedding
$\cone(Z)\to\prod_{a\in A}T_a$.
\qed

\section{Proof of Theorem~\ref{thm:main}}\label{sect:proofmain}

\subsection{Basics of hyperbolic spaces}\label{subsect:hypspaces}

We briefly recall necessary facts from the hyperbolic spaces
theory. For more details the reader may consult e.g. \cite{BoS}.

Let
$X$
be a metric space. Fix a base point
$o\in X$
and for
$x$, $x'\in X$
put
$(x|x')_o=\frac{1}{2}(|xo|+|x'o|-|xx'|)$.
The number
$(x|x')_o$
is nonnegative by the triangle inequality, and it is
called the Gromov product of
$x$, $x'$
w.r.t.
$o$.

A metric space
$X$
is {\em (Gromov) hyperbolic} if
the
$\de$-{\em inequality}
$$(x|x'')_o\ge\min\{(x|x')_o,(x'|x'')_o\}-\de$$
holds for some
$\de\ge 0$,
some base point
$o\in X$
and all
$x$, $x'$, $x''\in X$.

Let
$X$
be a hyperbolic space and
$o\in X$
be a base point. A sequence of points
$\{x_i\}\sub X$
{\em converges to infinity,} if
$$\lim_{i,j\to\infty}(x_i|x_j)_o=\infty.$$
Two sequences
$\{x_i\}$, $\{x_i'\}$
that converge to infinity are {\em equivalent} if
$$\lim_{i\to\infty}(x_i|x_i')_o=\infty.$$

The {\em boundary at infinity}
$\di X$
of
$X$
is defined as the set of equivalence classes
of sequences converging to infinity.
The Gromov product extends to
$X\cup\di X$
as follows. For points
$\xi$, $\xi'\in\di X$
the Gromov product is defined by
$$(\xi|\xi')_o=\inf\liminf_{i\to\infty}(x_i|x_i')_o,$$
where the infimum is taken over all sequences
$\{x_i\}\in\xi$, $\{x_i'\}\in\xi'$.
Note that
$(\xi|\xi')_o$
takes values in
$[0,\infty]$,
and that
$(\xi|\xi')_o=\infty$
if and only if
$\xi=\xi'$.

Similarly, the Gromov product
$$(x|\xi)_o=\inf\liminf_{i\to\infty}(x|x_i)_o$$
is defined for any
$x\in X$, $\xi\in\di X$,
where the infimum is taken over all sequences
$\{x_i\}\in\xi$.

A metric
$d$
on the boundary at infinity
$\di X$
of
$X$
is said to be {\em visual}, if there are
$o\in X$, $a>1$
and positive constants
$c_1$, $c_2$,
such that
$$c_1a^{-(\xi|\xi')_o}\le d(\xi,\xi')\le c_2a^{-(\xi|\xi')_o}$$
for all
$\xi$, $\xi'\in\di X$.
In this case, we say that
$d$
is the visual metric with respect to the base point
$o$
and the parameter
$a$.
The boundary at infinity is bounded and complete w.r.t.
any visual metric, and if
$a>1$
is sufficiently close to 1, then a visual metric with
respect to
$a$
does exist.

\subsection{Visual hyperbolic spaces}\label{subsect:visual}

A hyperbolic space
$X$
is called {\em visual}, if for some base point
$x_0\in X$
there is a positive constant
$D$
such that for every
$x\in X$
there is
$\xi\in\di X$
with
$|xx_0|\le(x|\xi)_{x_0}+D$
(one easily sees that this property is independent of
the choice of
$x_0$).

A map
$g:X\to Y$
to a metric space
$Y$
is said to be {\em roughly homothetic} if
$$|g(x)g(x')|\doteq\La|xx'|$$
up to a uniformly bounded additive error for some
constant
$\La>0$
and all
$x$, $x'\in X$.
The space
$X$
is said to be {\em roughly similar} to the image
$g(X)$.

For the proof of the following Proposition see
\cite[Proposition~6.2]{Bu}, it can also be extracted from \cite{BoS}.

\begin{pro}\label{pro:roughsim} Every visual hyperbolic space
$X$
is roughly similar to a subspace of the hyperbolic cone
over the boundary at infinity,
$\cone(\di X)$,
where
$\di X$
is taken with a visual metric.
\qed
\end{pro}

Combining Proposition~\ref{pro:roughsim} and
Theorem~\ref{thm:embcone}, we obtain a quasi-isometric
embedding
$X\to\prod_{a\in A}T_a$
with
$|A|=\cdim\di X+1$,
which completes the proof of Theorem~\ref{thm:main}.
\qed

%%%%%%%%%%%%%%%%%%%%%%%%%%%%%%%%%%%%%%%%%%%%%%%%%%%%%%%%%%%%%%

\bigskip
\begin{tabbing}

St. Petersburg Dept. of Steklov\hskip11em\relax \= \\

Math. Institute RAS, Fontanka 27, \> \\

191023 St. Petersburg, Russia\> \\

{\tt sbuyalo@pdmi.ras.ru}\> \\

\end{tabbing}

\end{document}